\documentclass[12pt]{article}
\usepackage{amssymb}

\newtheorem{thm}     {Theorem}[section]
\newtheorem{prop}    [thm]{Proposition}

\newtheorem{cor}     [thm]{Corollary}
\newtheorem{lemma}   [thm]{Lemma}

\newcommand{\proof} {\noindent{\bf Proof. }}
\newcommand{\A}{{\mathcal A}}
\newcommand{\B}{\mathbb B}
\newcommand{\C}{\mathbb C}
\newcommand{\D}{\mathbb D}
\newcommand{\Pro}{\mathbb P}
\newcommand{\R}{\mathbb R}

\newcommand{\st}{{\rm st}}

\def\Re{{\rm Re\,}}

\def\bar{\overline}

\begin{document}

\title{Gluing complex discs to Lagrangian manifolds by Gromov's  method  }
\author{Alexandre Sukhov{*} and Alexander Tumanov{**}}
\date{}
\maketitle

{\small

* Universit\'e des Sciences et Technologies de Lille, Laboratoire
Paul Painlev\'e,
U.F.R. de
Math\'e-matique, 59655 Villeneuve d'Ascq, Cedex, France, sukhov@math.univ-lille1.fr

** University of Illinois, Department of Mathematics
1409 West Green Street, Urbana, IL 61801, USA, tumanov@math.uiuc.edu
}
\bigskip

{\small Abstract. The paper discusses some aspects of Gromov's theory.

MSC: 32H02, 53C15.

Key words: almost complex structure, symplectic structure, Lagrangian manifold,
$J$-complex disc.

\bigskip

"Attachement des disques complexes aux vari\'et\'es Lagrangiennes par la m\'ethode de Gromov."

R\'esum\'e. L'article discute certains aspects de la th\'eorie de Gromov.

Mots-cl\'es: structure presque complexe,structure symplectique, vari\'et\'e Lagrangienne,
disque $J$-complex.}

\section{Introduction}

The present paper discusses some aspects of the  work of M.Gromov \cite{Gr} where the method of pseudo-holomorphic curves was introduced and successfully applied to several fundamental problems of symplectic topology; "...the most 
striking results in symplectic and contact topology have been so far obtained only by this method..." \cite{Ar}. The concentration of ideas in Gromov's paper is high and some of them are only sketched. Detailed proofs,  additional technical ingredients and far reaching generalizations have been elaborated by many authors enlarging an impressive area of applications. At present there exist several excellent introductions to the theory of pseudo-holomorphic curves focused on its different aspects and various applications,
see for instance  \cite{Aud,ElTr,IvSh,Hu,MS} (this list is, of course, highly incomplete even in the category of monographs and expository articles). A brief but deep description of Gromov's ideas is given in \cite{D}. Our modest goal  is  to present some of the  results of original Gromov's work mainly from the point of view of  Complex Analysis and PDE theory . This paper is not a survey so the references list is quite short; an interested reader can consult the above mentioned monographs  and the  exponentially  growing literature.

{\small Acknowledments. The first author is supported by Labex SEMPI. This work was partially done when the first author visited the Indiana University and the University of Illinois (Urbana-Champaign) during the Fall 2011 and the second author visiting Laboratoire Paul Painl\'ev\'e (USTL) during the Fall 2012 as CNRS researcher. We cordially thank these institutions for their support.

We are grateful to the referee for his suggestions improving the first version of the paper.}

\section{Preliminaries}

{\bf 1. Almost complex manifolds and their maps.} Let $M$ be a smooth ($C^\infty$) manifold of dimension $2n$. An almost complex structure $J$ on $M$ is a map (of appropriate regularity class $C^k$ or $C^\infty$) which associates to every point $p \in M$ a linear isomorphism $J(p): T_pM \to T_pM$ of the tangent space $T_pM$ satisfying $J(p)^2 = -I$,  $I$ being the identity map. A couple $(M,J)$ is called an almost complex manifold of complex dimension n.

Let $(M,J)$ and $(M',J')$ be smooth  almost complex manifolds. A $C^1$-map $f:M' \to M$ is called  $(J',J)$-complex or  $(J',J)$-holomorphic  if it satisfies {\it the Cauchy-Riemann equations} 
\begin{eqnarray}
\label{CRglobal}
df \circ J' = J \circ df.
\end{eqnarray}
 Gromov's theory is devoted to the case where $M'$ has the complex dimension 1; hence the structure $J'$ is necessarily integrable (see, for instance, \cite{Aud}) and $M'$ is a Riemann surface. In this special case 
holomorphic maps are called $J$-complex (or $J$-holomorphic) curves.
We use the notation  $\D$ for  the
unit disc in $\C$ and  $J_{\st}$ for  the standard complex structure
of $\C^n$; the value of $n$ will be clear from the context. If in the above definition we have $M' = \D$ and $J' = J_{st}$, we  call such a map $f$ a $J$-{\it complex  disc} or a
 pseudo-holomorphic disc or just a  holomorphic disc
if $J$ is fixed. Similarly, if $M'$ is the Riemann sphere, $f$ is called a $J$-{\it complex sphere}.

Let $(M,J)$ be an almost complex manifold and $E \subset M$ be a real submanifold of $M$. 
Suppose that a $J$-complex disc $f:\D \to M$ is  continuous on $\overline\D$  and 
satisfies $f(b\D) \subset E$. 
Then we say that (the boundary of ) the disc  $f$ is {\it glued} or {\it attached} to $E$ or simply 
that $f$ is attached to $E$. Sometimes such maps are called {\it Bishop discs} for $E$ and we employ this terminology. Of course, if $p$ is a point of $E$, 
then the constant map $f \equiv p$ always satisfies this definition. 
Often it is of interest is to prove an existence (or non-existence) of a {\it  non-constant} $J$-complex  disc attached to $E$. Gromov's theory provides a powerful tool for these studies.

\bigskip

{\bf 2. Cauchy-Riemann equations in coordinates.} In local coordinates $Z\in\C^n$, an almost complex structure
$J$ is represented by a $\R$-linear operator
$J(Z):\C^n\to\C^n$, $Z\in \C^n$ such that $J(Z)^2=-I$. We will use the notation $\zeta=\xi+i\eta\in\D$. Then the Cauchy-Riemann equations
(\ref{CRglobal}) for a $J$-complex  disc $Z:\D\to\C^n$, $Z: \D \ni \zeta \mapsto Z(\zeta)$ have the form $Z_\eta=J(Z)Z_\xi$. Similarly to  \cite{Aud},
we represent $J$ by a complex $n\times n$
matrix function $A=A(Z)$ so that the Cauchy-Riemann equations
have the form
\begin{eqnarray}
\label{holomorphy}
Z_{\bar\zeta}=A(Z)\bar Z_{\bar\zeta},\quad
\zeta\in\D.
\end{eqnarray}
We first discuss the relation between $J$ and $A$ for
fixed $Z$. Let $J:\C^n\to\C^n$ be an $\R$-linear map
so that $\det(J_\st+J)\ne0$.
Put $Q=(J_\st+J)^{-1}(J_\st-J)$. Then 
$J^2=-I$ if and only if
$QJ_\st+J_\st Q=0$, that is, $Q$ is complex anti-linear.

We introduce
$$
\mathcal J=\{J:\C^n\to\C^n: J\;{\rm is}\;\R{\rm-linear},\;J^2=-I,\;
\det(J_\st+J)\ne0  \}
$$
$$
\mathcal A=\{A\in Mat(n,\C): \det(I-A\bar A)\ne0 \}
$$

Let $J\in\mathcal J$. Then the defined above map $Q$  is anti-linear, hence, there is a unique
matrix $A \in Mat(n,\C)$ such that $Av=Q\bar v$, $v\in\C^n$. It is proved in \cite{Aud} (see also \cite{SuTu11a}) that the map $J\mapsto A$ is a birational
homeomorphism $\mathcal J\to\mathcal A$. We sum up. Let $J$ be an almost complex structure in a domain
$\Omega\subset\C^n$. Suppose $J(Z)\in\mathcal
J$, $Z\in\Omega$.
Then $J$ defines
a unique complex matrix function $A$ in $\Omega$ such that
$A(Z)\in\A$, $Z\in\Omega$. We call $A$ the
{\it complex matrix} of $J$ avoiding an employment the general Kodaira deformation theory terminology. The matrix $A$ has the
same regularity properties as $J$. Therefore, the notation $J_A$ or $A_J$ is appropriate according to the sense in which the correspondence $A \leftrightarrow J$ is viewed.

Let $M$ be an almost complex manifold of complex dimension $n$. Locally every almost complex structure $J$ on $M$ admits the complex matrix in a suitable coordinate chart. Denote by $\B_n$ the euclidean unit ball of $\C^n$.  For every point $p \in
M$, every $k \geq 1$  and every $\lambda_0 > 0$ there exist a neighborhood $U$ of $p$ and a
coordinate diffeomorphism $Z: U \to \B_n$ such that
\begin{eqnarray}
\label{firstnorm}
Z(p) = 0, \quad dZ(p) \circ J(p) \circ dZ^{-1}(0) = J_{st}
\end{eqnarray}  
and the
direct image $Z_*(J):= dZ \circ J \circ dZ^{-1}$ satisfies 
\begin{eqnarray}
\label{secondnorm}
\vert\vert Z_*(J) - J_{st}
\vert\vert_{C^k(\overline {\B_n})} \leq \lambda_0.
\end{eqnarray} Indeed, first consider 
 a diffeomorphism $Z$ between  a neighborhood $U'$ of
$p \in M$ and  $\B_n$  satisfying (\ref{firstnorm}). Then for $\lambda > 0$ introduce the isotropic dilation
$d_{\lambda}: t \mapsto \lambda^{-1}t$ in $\C^n$ and the composition
$Z_{\lambda} = d_{\lambda} \circ Z$. Clearly $ \vert\vert (Z_{\lambda})_{*}(J) - J_{st} \vert\vert_{C^k(\overline{\B_n})} \to 0$ as $\lambda \rightarrow
0$. Setting $U = Z^{-1}_{\lambda}(\B_n)$ for
$\lambda > 0$ small enough, we obtain a coordinate chart satisfying (\ref{firstnorm}), (\ref{secondnorm}). This elementary  observation  is often used in the local theory of $J$-complex curves.

\bigskip

{\bf 3. Analytic tools.}  The main analytic  tool in the theory of $J$-complex curves is the Cauchy-Green integral 
\begin{eqnarray}
\label{CGintegral}
T f(\zeta) = \frac{1}{2\pi i}\int_\D \frac{f(\tau)}{\tau - \zeta}d\tau \wedge d\overline\tau
\end{eqnarray}
 Denote by $C^{k,\alpha}(\D)$, $k \geq 0$, $0 < \alpha < 1$, the usual H\"older space of $C^{k,\alpha}$-functions in $\D$. A  classical property of $T$ is its regularity asserting that  $T:C^{k,\alpha}(\D) \to C^{k+1,\alpha}(\D)$ is a linear bounded operator. The  importance of the Cauchy-Green operator comes from the fundamental fact that $T$ gives a solution for the $\overline\partial$-equation in $\D$: we have $(Tf)_{\bar\zeta} = f$ for every $f \in C^{k,\alpha}(\D)$. As an example, consider 
the result of  Nijenhuis and Woolf (see for instance \cite{Aud}) which lies in the very foundation of theory. It 
states that for a given point $p \in M$ and a tangent
vector $v \in T_pM$ there exists a $J$-complex disc
$f: \D \to M$ such that $f(0)=p$ and
$df(0)(\frac{\partial}{\partial\xi}) = \lambda v$ for
some $\lambda> 0$. 
The disc $f$ can be chosen smoothly depending on the
initial data $(p,v)$ and the structure $J$. The above-mentioned regularity of $T$ allows to prove this theorem quite similarly to the Cauchy existence theorem for ODE's. Indeed, we replace the Cauchy-Riemann equations (\ref{holomorphy}) by the integral equation 
\begin{eqnarray}
\label{NW}
Z - T\left ( A(Z)\bar Z_{\bar\zeta} \right ) = W
\end{eqnarray}
where $W$ is a usual holomorphic vector-valued function in $\D$. One can assume that $A(0) = 0$ i.e. $J(0 ) = J_{st}$; as shown above, after isotropic dilations of coordinates the norm of $A$ is small.
But then the implicit function theorem establishes a one-to-one correspondence between the solutions $Z$  of the integral equation (\ref{NW})  and usual holomorphic discs $W$ in a sufficiently small neighborhood of the origin. This implies  the theorem. 

This simple principle is behind many local properties of $J$-complex curves. It  allows to develop their local theory explicitely, employing the classical properties of the Cauchy-Green integral and  related  singular integrals   without  the general elliptic PDE machinery.
This is due to the well-known particularity of  the theory of first order elliptic PDE with two independent variables and its interactions with  Complex Analysis (cf. \cite{BeS,Mo, Ve}).
As a consequence, the local properties of $J$-complex curves are similar to the properties of usual complex curves in complex manifolds (though complete proofs sometimes require substantial technical efforts). For example, the set of critical points of a non-constant $J$-complex curve is discrete and the intersection set of two  $J$-complex curves with distinct images is also discrete. Further important consequence is {\it the positivity of intersections property} and {\it the adjunction formula} for $J$-complex curves, see \cite{Aud,MS,MiWh}.

Using the generalized Cauchy formula and adding into the equation (\ref{NW}) suitable terms containing the usual Cauchy integral (over $b\D$), one can use such a modified integral equation  to obtain solutions to some boundary value problems for the Cauchy-Riemann equations (\ref{holomorphy}). This allows to construct $J$-complex discs with boundary data only locally. The global case is the subject of Gromov's theory and, as we will see, requires more advanced methods of non-linear analysis. Nevertheless, it is useful to keep in mind that in coordinates we are dealing with boundary value problems for the equations (\ref{holomorphy}).

\bigskip

{\bf 4. Interaction with symplectic and metric structures.} Let $M$ be a smooth real manifold of dimension $2n$. A closed non-degenerate exterior 2-form $\omega$ on $M$ is called a {\it symplectic form} on $M$. A couple $(M,\omega)$ is called a symplectic manifiold. As an example, consider $\C^{n}$ with the coordinates $z_j = x_j + iy_j$. The form $\omega_{st} = \sum_{j=1}^n dx_j \wedge dy_j = (i/2)\sum_{j=1}^n dz_j \wedge d\bar z_j$ is called {\it the standard symplectic form}. According to the classical Darboux's theorem \cite{Aud}, every symplectic form $\omega$ on $M$ is locally conjugated (or symplectomorphic) to $\omega_{st}$, i.e. there exists a local coordinate diffeomorphism $\phi$ satisfying $\phi^*\omega_{st} = \omega$. One of the consequence of Gromov's theory is that globally (on the whole $\R^{2n}$) this propety
fails.

Let $(M,J)$ be an almost complex manifold. A {\it $J$-hermitian metric} on $M$ is a real bilinear form $h:TM \times TM \to \C$ such that
\begin{itemize}
\item[(a)] $h(Ju,v) = i h(u,v) = i\overline h(v,u)$.
\item[(b)] $h(u,u) > 0$, $\forall u \neq 0$.
\end{itemize}

Given hermitian metric $h$, we have the decomposition on its real and imaginary parts: $h(u,v) = g_h(u,v) - i\omega_h(u,v)$. Then $g_h$ is a Riemannian metric on $M$ and $\omega_h$ is an exterior $2$-form. Furthermore, $g_h(u,v) = \omega_h(u,Jv)$
and $\omega_h(Ju,Jv) = \omega_h(u,v)$. In particular $h(u,v) = \omega_h(u,Jv) - i\omega_h(u,v)$. The form $\omega_h$ is called {\it the 2-form associated with $h$}.  An exterior 2-form on $(M,J)$ is called $J$-{\it calibrated} if 
\begin{itemize}
\item[(a)] $\omega(Ju,Jv) = \omega(u,v)$
\item[(b)] $\omega(u,Ju) > 0$, $\forall u \neq 0$
\end{itemize}
Each calibrated form defines a hermitian form if we set $h(u,v):= \omega(u,Jv) - i\omega(u,v)$. We say that a 2-form $\omega$ {\it tames} an almost complex structure $J$ if
$$\omega(u,Ju) > 0, \forall u \neq 0$$
i.e. only the above assumption (b) is imposed. It is known that
 calibrating (resp. tamed) almost complex structures on a given symplectic manifold form a non-empty contractible space \cite{Aud,Gr,MS}. A model example is provided by the standard symplectic form $\omega_{st}$ and the standard complex structure $J_{st}$ of $\C^n$. Clearly, $\omega_{st}$ is $J_{st}$-calibrated. We mention also a useful characterization of almost complex structures $J$ on $\C^n$ tamed by $\omega_{st}$ in terms of their complex matrices, see \cite{Aud}. Namely, $J$ is $\omega_{st}$-tamed if and only if for every $Z \in \C^n$ one has $\parallel A_J(Z) \parallel < 1$; here the operator norm is induced by the Euclidean inner product. Also, $J$ is calibrating if in addition its complex matrix $A_J$ is symmetric.

Suppose now that $J$ is tamed by $\omega$. Then we can define the Riemannian
metric $g(u,v) = \frac{1}{2}[ \omega(u,Jv) + \omega(v,Ju)]$. This metric is called {\it the canonical Riemannian metric} associated with $\omega$ and
$J$. In the case where $\omega$ is calibrated by $J$, it coincides with the
metric $\omega(\bullet,J\bullet)$.

Let $(M,\omega,J)$ be an almost complex  manifold with $J$-calibrated symplectic structure.  Let $X$ be a $J$-complex submanifold of $(M,\omega,J)$ (i.e. at every point  the tangent space of $X$ is $J$-invariant)  of complex  dimension $k$ with the canonical orientation. Denote respectively by $dV_X$ the volume form and by $vol_{2k}X$ the volume of $X$ induced by the canonical  metric $g$.  Then $dV_X = (1/k!)\omega^k\vert_X$; furthermore, if 
 $X$ is an oriented real $2k$-dimensional submanifold in $(M,\omega,J)$, then 
$(1/k!)\int_X \omega^k \leq vol_{2k}X$
and the equality (in the case of a finite volume) holds if and only if $X$ is $J$-complex.
This volume estimate  is called {\it the Wirtinger inequality}. 
As a consequence,  $J$-complex submanifolds are minimal and their volume with respect to the metric $g$ is given by
$$vol_{2k} X = \frac{1}{k!}\int_X \omega^k$$

 Let $(M,\omega,J)$ be a tamed  almost complex manifold. Consider  a Riemann surface $(S,J_S)$ 
and  a $J$-complex curve $f: (S,J_S) \to (M,\omega,J)$. Its $\omega$-area (or {\it symplectic area}) is defined by
\begin{eqnarray}
\label{area}
area(f) = \int_S f^*\omega
\end{eqnarray}
If additionally $J$ calibrates $\omega$, this is precisely the area defined by the metric $g$ associated with $\omega$ and $J$ and  $f(S)$ is a minimal surface for this metric. In the tamed case  the minimality in  general fails. Consider the special case where $S = \D$ i.e. $f$ is a $J$-complex disc in a  tamed almost complex manifold. Then
the expression 
$$E(f):=  \frac{1}{2}\int_\D \left ( \left \vert \left \vert \frac{\partial f}{\partial
  \xi} \right\vert \right \vert^2_g + \left\vert\left\vert \frac{\partial f}{\partial
  \eta} \right\vert\right\vert^2_g \right )d\xi \wedge d\eta$$
where the norm $\parallel \bullet \parallel_g$ is taken with respect $g$, is called {\it the energy} of $f$. We have
$$E(f) = \int_\D f^*\omega$$
This  equality is called the energy identity, see for instance \cite{MS}.

{\bf Example.} As a consequence  every non-constant $J$-complex curve has a strictly positive symplectic  area. For instance, suppose that $\omega$ is globally exact on $M$. Then for every almost complex structure tamed by $\omega$ the manifold $(M,\omega,J)$ does not contain non-constant $J$-complex spheres. Indeed, by Stokes' formula the symplectic area of such a sphere is equal to  $0$.

\bigskip

{\bf 5. Real submanifolds.} Let $(M,\omega,J)$ be a tamed manifold of complex dimension $n$. A real submanifold $E$ of real dimension $n$ in $M$ is called

\begin{itemize}
\item[(i)] {\it Lagrangian} if $\omega\vert_E = 0$.
\item[(ii)] {\it Totally real} if $T_pE \cap J(p)(T_pE) = \{ 0 \}$ for every $p \in E$.
\end{itemize} 
It is well-known that every Lagrangian submanifold is totally real as well as that the inverse in general fails. For example, the standard torus $\Lambda = b\D \times ...\times b\D = (b\D)^n$ in 
$(\C^n,\omega_{st},J)$ is Lagrangian and  totally real for every $J$ tamed by $\omega_{st}$. This example will be in the focus of our study.

Let $E$ be a Lagrangian or totally real submanifold in $M$. Suppose that  $E$ is  the zero set $E = \rho^{-1}(0)$ of a smooth vector function $\rho: M \to \R^n$. If $f$ is a Bishop disc for $E$, then 
\begin{eqnarray}
\label{boundary}
\rho \circ f(\zeta) = 0, \zeta \in b\D
\end{eqnarray}
In local coordinates $f$ satisfies the Cauchy-Riemann equations (\ref{holomorphy}) which together with (\ref{boundary}) form a non-linear  elliptic boundary value problem. An appropriate tool of   the non-linear analysis here is the continuity method. The strategy is the following.  Given boundary value problem one associates   a homotopy in the suitably choosen spaces of PDE operators (i.e. essentially the complex matrices $A$ in our case) and the boundary value data (i.e. the above  functions $\rho$) joining the initial problem with a simpler one for which a solution can be constructed. Next one constructs a homotopy in the space of solutions in order to go back to the initial problem and to obtain its solution. This procedure is based on  two main technical ingredients. 

\bigskip

The first one is an analysis of the linearized boundary value problem which allows to  extend slightly by the implicit function theorem  the homotopy path in the space of solutions. Usually the Fredholm properties of the linearized problem are useful here. In our case they again follow essentially from the regularity properties of the Cauchy-Green integral. This is one of the central result of the classical theory of linear singular integral equations developed  from
50-s to 70-s (first for a scalar equation on the plane, and then for vector-valued dependent variables). The book \cite{Wend} contains a rather complete survey of this theory. For reader's convenience we include to Section 7 a short proof of the Fredholm property for the model boundary value problem with the usual $\overline\partial$-operator. Substantially more general operators are considered in \cite{Wend}. An application of standard methods based on the Cauchy integral theory requires a global coordinate neighborhood for a prescribed $J$-complex disc. If such a disc  is not embedded or immersed, one can consider its graph and a suitable lift of an almost complex structure. In \cite{SuTu12} this approach is used in order to consruct the deformation theory for $J$-complex discs with free boundaries. However, a study of the Bishop discs requires a deformation theory with Lagrangian or totally real boundary data. In the case of complex dimension 2 this is rather simple. Given $J$-complex disc glued to a Lagrangian or totally real manifold, one can associate an integer invariant under homotopy: {\it the Maslov index}, see for instance \cite{MS}. An existence of nearby discs (under a perturbation of an almost complex structure and  boundary data), as well as the maximal number of real variables parametrizing the perturbed discs is completely determined by this index, see for instance \cite{GaSu,MS,HLS}. Essentially this is a direct consequence of the classical theory of   the  linear Riemann-Hilbert boundary value problem in the unit disc for usual (or generalized) scalar analytic  functions. The number of parameters in the general solution is  completely determined by the winding number of the coefficient from the boundary value condition, see \cite{Ve}.
Another approach is  especially fruitful in the case of  compact $J$-complex curves.  One considers the pull-back  of the tangent bundle of $M$ by a given $J$-complex map. This gives rise to a complex vector bundle over the source Riemann surface and the Cauchy-Riemann equations can be expressed intrinsically in terms of associated  connections and metric structures. In this way the well-elaborated machinery of elliptic operators on vector bundles can be applied. This general approach is employed by many authors, see \cite{IvSh,Gr,MS}. In the case of complex bundle of rang 2 (corresponding to the complex dimension 2 of the target almost complex manifold) a deformation of a given compact $J$-complex curve is again determined by a single homotopy invariant: the first Chern class of the bundle. Unfortunately, both in the compact or Lagrangian boundary data case the situation changes seriously when the complex dimension of $M$ is higher than $2$. Though the first Chern class or,respectively, the Maslov index are still defined, a possibility of deformation of a given $J$-complex curve depends on a finite number of additional charactersitics which in general are not stable under homotopy. For instance, in the model case of the linear Riemann-Hilbert boundary value problem for usual vector-valued analytic functions the solvability depends on the  so called partial indices which are not homotopically stable, see for instance \cite{MiPr}. A similar problem occurs in the compact case. This difficulty was overcomed in Gromov's theory by geometrization of the problem using  the Sard-Smale theorem. This explains the substantial difference between the theory in  complex dimension 2 and higher dimensions.

\bigskip

The second ingredient is {\ a priori}  estimates often coming from geometric considerations. In our case they are incorporated into   Gromov's compactness theorem giving a very strong convergence of sequences of $J$-complex curves with uniformly bounded areas.

\bigskip

{\bf 6. Gromov's Compactness theorem.}  Let $S$ be a compact Riemann surface with (possibly empty) smooth boundary $bS$.
We use the canonical identification of the
complex plane $\C$ with $\C\Pro \backslash \{ \infty \}$. Let
$(M,\omega,J)$ be a symplectic manifold with a tamed almost complex
structure; as above, $g$ is the associated Riemann metric. We assume that $M$ has bounded geometry with respect to $g$ i.e. satisfies the standard assumptions on the completeness, curvature and injectivity radius, see \cite{Aud}, p.178. Let $E$ be a smooth compact totally
real submanifold of maximal dimension in $M$.

Consider a sequence $f^n: S \to M$ of $J$-complex  maps such that $f^n(bS) \subset E$.

Let $\psi: \C\Pro \to M$ be a {\it non-constant} $J$-complex map. We
say that $\psi$ occurs as a {\it spherical bubble} for the sequence $(f^n)$ if there exists a sequence of holomorphic charts $\phi^n :R_n\D \to S$ with $R_n \to \infty$
converging uniformly on compacts subsets of $\C$ to a point $p \in S$ and such 
that   $f^n \circ \phi^n \to \psi$ uniformly on compact subsets of $\C$.

Let $\psi:\D \to M$ be a {\it non-constant} $J$-complex map, continuous on
$\overline\D$, with $\psi(b\D) \subset E$. We say that $\psi$ occurs as a {\it disc bubble} for the sequence $(f^n)$ if
there exists a sequence of holomorphic charts $\phi^n:\D \backslash (-1 + \delta_n\D) \to
  S \cup bS$, smooth on $\overline\D \backslash ( -1 + \delta_n\D )$ with 
$\phi^n (b\D   \backslash ( -1 + \delta_n\D )) \subset E$ and
$\delta_n \to 0$, such that $(\phi^n)$ converge uniformly on compact subsets
of $\overline \D \backslash \{ - 1 \}$ to a point $p \in bS$ and $f^n \circ \phi^n \to \psi$ uniformly on compact subsets of $\overline \D \backslash \{ - 1 \}$.

One of the simplest versions of Gromov's Compactness
Theorem is the following:

\begin{prop}
\label{Gromov1}
Let $(f^k):S \to M$ be a sequence of $J$-complex maps continuous on $S
\cup bS$, $f^k(b S) \subset E$, intersecting a fixed compact
subset $ K \in M$ and such that 
$$area (f^k) \leq C$$
where $C > 0$ is a constant. Then there exists a finite set $\Sigma$ in $S \cup bS$, possibly empty, such that after extraction a subsequence we have:
\begin{itemize}
\item[(i)] $(f^k)$ converges uniformly on compact subsets of $(S \cup b S) \backslash \Sigma$ to a
  $J$-complex map $f^\infty:S \to M$. Furthermore, the convergence is in every $C^r$-norm.
\item[(ii)] A spherical bubble occurs at every point in $\Sigma \cap S$,
\item[(iii)] A  disc occurs at every point in $\Sigma \cap bS$.
\end{itemize}
\end{prop}

This result is sufficient for our goals. Much more advance versions and detailed proofs are contained in \cite{Aud,IvSh1,Hu,MS}. We conclude by some remarks.

{\bf 1.} It follows from the above definition of a spherical or disc bubble 
that if they arise, then  they are {\it  non-constant}. Therefore, if for 
some topological reasons non-constant bubbles can not arise, a sequence of $J$-complex maps uniformly converges in the closed unit disc. This  is often used in order to prove a convergence of sequences of $J$-complex curves.

{\bf 2.} "Preservation of energy". As above, $\Sigma = \{ p_1,...,p_l \}$ denotes the set of points where  bubbles arise. Given $\varepsilon > 0$ and $p_j \in \Sigma$ set 
$$m_\varepsilon(p_j) = \lim_{k\to\infty}E(f^k\vert_{p_j + \varepsilon\B_n})$$
and 
$$m(p_j) = \lim_{\varepsilon \to 0} m_\varepsilon(p_j)$$
Then 
$$E(f^\infty) + \sum_{j=1}^lm(p_j) = \lim_{k\to\infty} E(f^k)$$

{\bf 3.} Consider the special case $S = \D$. The sequence of sets $f^k(\bar \D)$ converges to a finite connected union $X$ of $J$-complex spheres and discs (a cusp-curve) in the Hausdorff distance. The above definition of bubbles depend on the initial parametrization of $f^n$ i.e. the above version of Gromov's theorem describes the convergence of maps, not the sets. For instance, after a suitable reparametrization of $S$ by a sequence of conformal isomorphisms, we can obtain a sequence $(\tilde{f}^n)$ converging outside the finite bubbling set to a single point in $E$. Then all limit spheres and discs in $X$ will occur as bubbles for  $(\tilde{f}^n)$.

{\bf 4.} The above theorem still holds, of course, if we vary the structure $J$ together with maps, i.e. $J$ is the limit in an appropiate $C^s$-norm of the sequence $(J^n)$ of almost complex structure and every $f^n$ is a $J^n$-complex map.

The proof of the above theorem is based on two  elementary tricks: the covering argument due to Sacks-Uhlenbeck and the renormalization argument essentially appearing in the definition of a bubble. This allows to show an arising of a finite number of finite area bubbles  defined on the punctured Riemann sphere or on the disc with punctured boundary respectively. In order to remove these isolated singularities, the standard elliptic estimates and the bootspraping argument can be applied. In the case of the disc bubble, these arguments are related to a suitable non-analytic version of the reflection principle.

\section{Gromov-Hartogs Lemma in complex  dimension 2}

In this section we solve the model boundary value problem for $J$-complex discs attached to a Lagrangian torus in $\C^2$. The small dimension allows to  control effectively the geometric properties of solutions.

We use the notation $Z = (z,w) \in \C^2 = \C \times \C$
for the standard coordinates in $\C^2$.  Set $\omega_1 =\frac{i}{2} d z \wedge d \bar{z}$ and $\omega_2 =   \frac{i}{2}d w \wedge d \bar{w}$ and  denote by $\omega = \omega_1 +  \omega_2$
the standard symplectic form on $\C^2$.

Let $M_2$ denotes the vector space of    complex
$(2 \times 2)$-matrix functions  $A$
defined on $\C^{2}$ and of class $C^{\infty}( \C^{2})$. Consider  the maps $A \in M_2$ satisfying  the following assumptions:

\bigskip

(i) The {\it strong  taming assumption} consists of two parts. First, we suppose that there exists a real  $0 \leq a_0 < 1$ such that

\begin{eqnarray}
\label{norm0}
\parallel A(Z) \parallel  < a_0, \forall Z \in \C^{2}
\end{eqnarray}
where the matrix norm is induced by the Euclidean norm of $\R^{4}$. Second, we suppose that 
{\it the map $Z \mapsto A(Z)$ is uniformly continuous on $\C^2$.}

\bigskip

Recall that the defined in the previous section map $J_A \leftrightarrow A$  establishes a one-to-one correspondence   between   matrix functions $A$ satisfying (\ref{norm0}) and  almost complex structures $J_A$ on $\C^{2}$  tamed by the standard symplectic form $\omega$. In what follows we simply denote $J_A$ by $J$. Our assumptions on $A$ and the explicit formula expressing $J_A$ in terms of $A$, see \cite{SuTu11a}, imply  that $J$ is uniformly continuous on $\C^2$. This guarantees that $(M,\omega,J)$ has the bounded geometry and allows to employ Gromov's compactness theorem. 

 The second assumption is 

\bigskip

(ii) The map $\C \ni z \mapsto (z,0)$ is $J$-complex. Writing explicitely

\begin{eqnarray}
\label{A-abc}
A= \left(
\begin{array}{cll}
a & & d\\
b & & c
\end{array}
\right)
\end{eqnarray}
we see that this assumption is equivalent to the condition $a(z,0) = b(z,0) = 0$ for every $z \in \C$.

\bigskip

Introduce the real 2-torus $\Lambda^t=b\D \times tb\D$ where $t > 0$. Then every torus $\Lambda^t$ is Lagrangian with respect to $\omega$. Therefore, $\Lambda^t$ is totally real with respect to each almost complex structure tamed by $\omega$.

\begin{thm}
\label{theo2}
Fix $t = T$. Under the above assumptions (i),(ii) the following holds.
\begin{itemize}
\item[(a)]
For every point $p = (1,q) \in\Lambda^T$ there exists a  $J$-complex disc
$f:\D\to \C^2$ of class $C^{\infty}(\bar\D)$ 
such that $f(1)=p$, $f$ is an embedding, $f(b\D)\subset \Lambda^T$,
and $f(\bar\D)$ does not meet $\bar\D\times\{0\}$.
Furthermore, $area(f) = \pi$. 
\item[(b)] When $q$ runs over the unit circle, the discs in (a) form a $C^{\infty}$-smooth one-parameter family. They are disjoint and fill a smooth Levi-flat (with respect to $J$) 
hypersurface $\Gamma \subset \C^2$
with boundary $\Lambda^T$. Furhermore, they depend continuously on $J$ and $t$.
\end{itemize}
\end{thm}

Since the proof is  short, we directly present it. Then we discuss relations of Theorem \ref{theo2} with other results.

\subsection{Proof of Theorem \ref{theo2}.}  The proof is based on the continuity method  discussed above. We proceed in several steps. Without loss of generality  assume $T = 1$ and write $\Lambda = b\D \times b\D$.

\bigskip

{\bf (1)} Since the torus $\Lambda$ is fixed, with some abuse of notation we denote again by $t$ the parameter which determines a homotopy of almost complex structures. Namely, for $t \in [0,1]$ consider the matrix $tA$ and the corresponding almost complex structure $J_t:= J_{tA}$. As a consequence of  the assumption (ii) the line $\C \times \{ 0 \}$ remains $J_t$-complex for all $t$.  Next, $J_0 = J_{st}$ and for $t = 0$ we have the Levi-flat hypersurface $\D \times b\D$  foliated by the embedded $J_{st}$-complex discs of the form  $h_c: \D \ni \zeta \mapsto (\zeta,c)$ where $c \in b\D$ is a constant.This provides for $t=0$ the discs $f$ with required properties. Suppose that the family of $J_t$-complex discs is defined on  $[0,t_0[$ with $0 \leq t_0 < 1$. Our goal is to extend this family with respect to the parameter $t$ on the whole interval $[0,1]$. 
\bigskip

{\bf (2)} Using the notation $z = x + iy$ and $w = u + i v$,   set $\lambda_1=  (1/2)(xdy - ydx)$, $\lambda_2 = (1/2)(udv - vdu)$ and $\lambda = \lambda_1 + \lambda_2$. Hence $\omega = d\lambda$. By continuity in $t$,   the restrictions of $z$- and $w$- components of the constructed above discs $f:b\D \ni \zeta \mapsto  f(\zeta)= (z(\zeta),w(\zeta))$ have the winding numbers about the origin equal to 1 and 0 respectively.  Since the components of $f$ take $b\D$ to the circles around the origin, by Stokes' formula we obtain $$area(f) = \int_{b\D} f^*\lambda = \pi.$$

Furhermore, the $z$-component vanishes somewhere in $\D$ since the winding number is equal to $1$. Reparametrizing the disc $f$ by a conformal automorphism of $\D$, we can assume that $z(0) = 0$ for all $t$. Next choose a point $p = (1,q) \in \Lambda$; one can also  assume  that $f(1) = p$. These {\it normalization conditions} define uniquely a family of $J_t$-complex discs $(f_t)$ depending continuously on $t \in [0,t_0[$. 

\bigskip

{\bf (3)} The key argument is the following 

\begin{lemma}
Suppose that every  disc $f_t: \zeta \mapsto (z_t(\zeta),w_t(\zeta))$, $t \in [0,t_0[$ does not intersect the axis $\C \times \{ 0 \}$.Then  there exists $\eta > 0$ such that 
$$\vert w_t(\zeta) \vert \geq \eta, \forall \zeta \in \bar\D, \forall t \in [0,t_0[$$
\end{lemma}
\proof  Suppose by absurd that there exists a sequence $f^k = (z^k,w^k)$, $f^k = f^{t(k)}$, $t(k) \to t_0$, such that $\inf_\D \vert w^k \vert \to 0$ as $k \to \infty$. The area of all discs are equal to $\pi$, the structures $J_t$ are tamed by assumption (i)  and the torus $\Lambda$ is totally real. Gromov's compactness theorem  implies (after extracting a subsequence) that the images $f^k(\D)$  converge in the Hausdorff metric to a finite  union of $J_{t_0}$-complex discs with boundaries  glued to the torus $\Lambda^{t_0}$. Notice here that $\omega$ is globally exact on $\C^2$ so every $J$-complex sphere in $(\C^2,\omega,J)$ is constant. Hence, spherical bubbles do not occur. Given such a limit disc, after a suitable reparametrization by a sequence of conformal isomorphisms, the convergence is in every $C^l(K)$-norm on each compact subset $K$ of $\bar\D \setminus \Sigma$ where $\Sigma$ is a finite subset of $b\D$. By assumption, one of the limit discs touches the $J_{t_0}$-complex line $\C \times \{ 0 \}$. Since the boundaries of discs is attached to $\Lambda^{t_0}$, the 
disc is not contained in this line. Positivity and stability of the intersection indices  of $J$-complex curves imply that $f^k(\D)$ also intersects the line $\C \times \{ 0 \}$ for $ k$ big enough. This  contradiction proves the lemma. $\blacksquare$

\bigskip

{\bf (4)} As a consequence we obtain that  the above sequence $f^k = (z^k,w^k)$ converges to a single $J_{t_0}$-complex disc in every $C^l$-norm on $\bar\D$. Indeed, consider  the "principal" limit disc $f^\infty = (z^\infty,w^\infty)$ of this sequence  defined as the  limit of the sequence of maps $f^k$ converging on $\bar\D$ (without any additional reparametrization) off at most a  finite subset of $b\D$ where bubbles arise. The disc $f^\infty$ has a non-constant $z$-component because of the normalization condition imposed above. In particular, its  area is positive. Since the winding number of its $w$-component is equal to zero by Lemma, we conclude that $area(f^\infty) = \pi$.  But  the total area of limit discs  is bounded $\pi$. Therefore  the limit set consists of a single disc and there  are no bubbles. Gromov's compactness theorem implies  the convergence  in every $C^l(\bar\D)$-norm. In particular, the winding numbers of the $z$- and $w$- components of the limit disc $f^\infty$ again are equal to $1$ and $0$ respectively. The discs under consideration are embeddings for $t \in [0,t_0[$.
Suppose that the limit disc $f^\infty$   is  multiply covered  that is $f^\infty = \tilde f \circ \Pi$ where $\tilde f$ is a $J_{t_0}$-complex disc, $\tilde f(b\D) \subset \Lambda$, and $\Pi$ is the Blaschke product of degree $d \geq 2$. Then the winding number of the $z$-component of $f^\infty$ is an integer multiple of $d$ and cannot be equal to 1. 
Thus the disc $f^\infty$ is not multiply covered and remains  an embedding  by the adjunction formula for $J$-complex curves. Then constructed discs remain embeddings for all $t$. Our families of discs and  almost complex  structures  are homotopic to the above  $J_{st}$-complex disc $h_c$ glued to the standard tori and  by continuity the Maslov index of every disc has the same value as for $h_c$ and so is equal to $0$. By the implicit function  theorem (see, for instance \cite{GaSu,HLS}) the disc $f_{t_0}$ generates a real 1-parameter family of  $J_t$-complex  discs with boundaries glued to $\Lambda$ for $ t \in [0,t_0 + \varepsilon[$ for some $\varepsilon > 0$ and satisfying the normalization condition (the discs $f_t$ already defined for  $t < t_0$ belong to this generated family by the uniqueness part of the implicit function theorem). This proves the part (a) of Theorem \ref{theo2}.
\bigskip

{\bf (5)} Consider another point $p' \in \Lambda^1$ and corresponding family of discs constructed as above. The discs of the two constructed  families do not intersect for $t$ close to $0$ and hence for all $t$ because of the positivity and stability of intersection indices of $J$-complex curves. This  implies (b) and concludes the proof of Theorem \ref{theo2}. $\blacksquare$

\subsection{Comments and remarks.}

{\bf 1.} J.Duval and D.Gayet \cite{DuGa} recently constructed an example of a totally real 2-torus in the unit sphere $S^3$ of $\C^2$, isotopic to the standard torus (i.e. unknotted) and such that there does not exist a $J_{st}$-complex disc with boundary attached to this torus. They proved that to every torus of this class one can attach the boundary of $J_{st}$-complex disc or the boundary of a $J_{st}$-complex annulus.

\bigskip

{\bf 2.} Theorem \ref{theo2} implies Gromov's non-squeezing theorem \cite{Gr} in $\C^2$. It suffices to use the discs provided by Theorem \ref{theo2} instead of  $J$-complex spheres in  classical Gromov's argument. We give the details in the next section.

\bigskip

{\bf 3.} Theorem \ref{theo2} can be used for studying holomorphically convex hulls (of course, this question is of interest only if the structure $J$ is integrable). For such application it is appropriate to give a slightly different construction.

 We  replace (ii) by a stronger assumption:

(ii') There exists $r_0 > 0$ such that  $A$ is a lower triangular complex matrix function i.e.
$d = 0$ in (\ref{A-abc}) on $\D \times r_0\D$ and   $\parallel a \parallel_\infty  \leq a_0$, $\parallel c \parallel_\infty \leq a_0$ for some $a_0 < 1$. Furthermore, we still assume that the map $\C \ni z \mapsto (z,0)$ is $J$-complex. Consider the tori $\Lambda^t = b\D \times t b\D$. We will construct a homotopy of $J_A$-complex discs starting from small $t$ and then extend it for all $t \in [0,1]$  without deformation of the almost complex structure $J_A$.

Reparametrizing our homotopy of tori in $t$ if necessary, we can  choose $r_1 > 0$  very small with respect to $r_0$ and  assume that $\Lambda^0 =   b\D \times r_1 b\D$ and $T = 1$. Then in view of the assumption (ii') it follows  by \cite{SuTu1} that there exists an open  neighborhood $W$ of $0$ in $\C$ such that 
the set $(\overline \D \times \bar W) \setminus (\overline \D \times \{ 0 \})$ is foliated by smooth Levi-flat (with respect to $J_A$) hypersurfaces whose boundaries coincide with the tori $\Lambda^t$, $t \in [0,r_2[$, $r_1 < r_2 < r_0$.  Every hypersurface in turn is foliated by $J$-complex discs $f$ homotopic to the  $J_{st}$-complex disc $h_c$ as above. Now we extend the homotopy in $t$ precisely as in the proof of Theorem.

\bigskip

{\bf 4.}  Denote by $H$ the unbounded "Hartogs type" figure $H = (\{ 1 - \delta < \vert z \vert < 1 \}  \times \C) \cup ( \D \times \varepsilon \D)$ with $0 < \delta < 1$, $\varepsilon > 0$. Assume that the Levi form (with respect to $J_A$) of the boundary $\Pi = b\D \times \C$ of the domain  $\D \times \C$ is non-negative definite  at every point. As a consequence of Remark 3 we obtain that this domain coincides with the holomorphic hull of $H$. We point out that the Levi flat hypersurfaces constructed in Theorem \ref{theo2} cannot touch the hypersurface $\Pi$ because of the Levi non-negative definiteness of $\Pi$ (see \cite{DiSu}). Furthermore,since $J$ is uniformly continuous by (i) and the area of discs are bounded, there exists an upper bound on diameters of discs with boundaries glued to a fixed torus (see, for instance, \cite{Aud}). This implies that the constructed Levi flat hypersurfaces sweep out the domain $\D \times \C$ i.e. the holomorphic envelope  of $H$. For this reason we call Theorem \ref{theo2}
"Gromov-Hartogs lemma".

\bigskip

{\bf 5.} M.Gromov \cite{Gr} proved (together with many other things) a result similar to Theorem \ref{theo2} in a more general setting of elliptic structures (an almost complex structure can be viewed as a special case). He assumes that a Lagrangian torus is contained in a strictly pseudoconvex hypersurface (this notion can be defined in the elliptic category). It is not obvious how to deform such a hypersurface keeping the strict pseudoconvexity through a deformation of almost complex (or elliptic) structure. However, if we apply the construction described in Remark 4, we do not need to deform a complex structure. This leads to the following bounded version  of the above  result. Suppose additionally that under the assumptions of  Remark 3 the torus $\Lambda = b\D \times b\D$ is contained in a smooth hypersurface which bounds a domain $\Omega$ and whose  Levi form with respect to $J_A$ is non-negative definite. Then the torus $\Lambda$ bounds the Levi-flat hypersurface (produced by the above construction) contained in $\Omega$. In this case it suffices to require that $A$ and $J$ are defined only in a neighborhood of the closure $\bar\Omega$.

\section{Gromov's Non-Squeezing Theorem} As above, let $\omega$ denote the standard symplectic form of $\C^2$. 

\begin{thm}
\label{squiz}
Let $G$ be a relatively compact domain  in $R\D \times \C$ where $R > 0$. Suppose that $r> 0$ and there exists a diffeomorphism $\Phi: r\B \to G$ with $\Phi^* \omega = \omega$. Then $r \leq R$.
\end{thm}
\proof Performing a translation in the $w$-direction one can assume that the disc $R\D \times \{ 0 \}$ does not meet $\overline G$. Consider an increasing  sequence $r_n \to r$. The almost complex structure $J := \Phi_* (J_{st})$ is tamed by $\omega$. Multiplying the complex matrix of $J$ by  suitable smooth cut-off functions, we obtain for every $n$ a smooth almost complex structure $J_n$ on $\C^2$ such that $J_n = J$ on $\Phi(r_n\B_2)$ and $J_n = J_{st}$ on $\C^2 \setminus \overline G$. Consider the point $p = \Phi(0)$. According to previous section, for every $n$ there exists a $J_n$-complex disc $f^n$ such that 

\begin{itemize}
\item[(i)] $f^n(0) = p$, 
\item[(ii)] $f^n(b\D) \subset Rb\D \times tb\D$ for some $t > 0$,
\item[(iii)]$area(f^n) = \pi R^2$
\end{itemize}

Then $X^n = \Phi^{-1}(f^n(\D) \cap \Phi(r_n\B_2))$ is a closed $J_{st}$-complex curve  in $r_n\B_2$. Furthermore, $0 \in X^n$ and $area(X^n) \leq \pi R^2$. Passing to the limit as $n \to \infty$, we obtain by Bishop's convergence theorem \cite{Ch} that there exists a closed complex curve $X$  in $r\B_2$, containing the origin and with  $area(X) \leq \pi R^2$.
On the other hand, since $0 \in X$, it follows by the classical results \cite{Ch} that $area(X) \geq \pi r^2$ and theorem follows.

\section{Gromov-Hartogs Lemma in higher  dimension}

An attempt to generalize directly the previous argument to higher dimensions meets difficulties. For instance, the positivity of intersections does not make sense. Another problem concerns a possibility of deformation of a $J$-complex disc with a Lagrangian (or totally real) boundary data. As it was discussed above,  such a deformation can not be described in terms of a single index invariant under a homotopy. So more advanced tools are needed. 

\bigskip

We use the notation $Z = (z,w) = (z,w_2,...,w_{n}) \in \C \times \C^{n-1}$
for the standard coordinates in $\C^n$.  Set $\omega_1 =\frac{i}{2} d z \wedge d \bar{z}$ and $\omega_j =   \frac{i}{2}d w_j \wedge d \bar{w}_j$. Let $\omega = \sum_{j=1}^n \omega_j$ denotes 
the standard symplectic form on $\C^n$. Let $G$ be an open set in $\R^N$, let also  $0 < \alpha  < 1$ and  $ k\geq 0$ be
an integer. We denote by $C^{k,\alpha}(G)$ 
the class of  functions $u:G \to \R$ admitting the partial derivatives $D^su$, $\vert s \vert \leq k$ in $G$ which are  $\alpha$- Holder continuous on $G$ when $\vert s \vert = k$. This is a Banach space with respect to the standard norm. For simplicity of notations we keep the same notation for the space of vector-valued functions $u:G \to \R^{N'}$  with components of class $C^{k,\alpha}(G)$.

 Denote by $M_n$ the vector space of    complex
$(n \times n)$-matrix functions  $A$
defined on $\C^{n}$ and of class $C^{k,\alpha}( \C^{n})$. Consider  the maps $A \in M_n$ satisfying  the following assumptions:

\bigskip

(i) The {\it strong taming assumption} is precisely the same as in dimension 2.

\bigskip

The next condition we impose  is

(ii) {\it (Support assumption) }  The support $K:= supp A$ is  separated from the union of hyperplanes  $\{ (z,w) \in \C^n: w_j = 0 \}$, $j=2,...,n$.  Therefore, the almost complex structure $J_A$ coincides with the standard complex structure $J_{st}$ of $\C^n$ in a neighborhood of this union. More precisely, there exists  constant $r_0 > 0$  such that the following separation property holds: 
\begin{eqnarray}
\label{support}
\inf \{ \vert w_j \vert : (z,w) \in K, j=2,...,n \} \geq 2r_0 
\end{eqnarray}
We point out that $supp A$ in general is not supposed to be a compact subset  in $\C^{n}$. For example,   it can be unbounded and  we do not require that   the structure $J_A$ could be extended to an almost complex structure on the complex projective space.

\bigskip

 We study  maps $Z :\zeta \mapsto Z(\zeta) = (z(\zeta),w(\zeta))$, $Z: \D \to
\C^n$ of  class $C^{k+1,\alpha}(\D)$, $k \geq 0$, satisfying  the following elliptic PDE system:

\begin{eqnarray}
\label{CR1}
Z_{\bar\zeta} - A(Z)\bar Z_{\bar\zeta} = 0, \zeta \in \D
\end{eqnarray}
with the non-linear boundary value condition
\begin{eqnarray}
\label{BdCn}
\vert z(\zeta) \vert^2   = R, \vert w_j(\zeta) \vert^2  = t_j, j = 2,...,n, \zeta \in b\D
\end{eqnarray}
Here $R > 0$, $t_j > 0$ are prescribed constants and the matrix function $A \in C^{k,\alpha}(\C^{n})$ satisfies the assumptions (i), (ii). A $C^1$-map $Z:\D \to \C^n$ is a $J_A$-complex disc if and only if it satisfies the system (\ref{CR1}) which represents the Cauchy-Riemann equations corresponding to the structure $J_A$.

The main result of this section is the following :

\begin{thm}
\label{theo1}
For every point $(a,b) \in \D \times \C^{n-1}$, $b = (b_2,...,b_n)$, $b_j \neq 0$, there exists $t \in
(R_{+}^*)^{n-1}$ and a solution $Z = (z,w) \in C^{k+1,\alpha}(\D)$ of the boundary value problem (\ref{CR1}), (\ref{BdCn}) satisfying $Z(0) = (a,b)$. Furthermore, the winding
number of the $z$-component of $Z$  is equal to $1$ and  the winding number of every $w_j$-component  is equal
to $0$.
\end{thm}

This result has several applications. We mention some of them.

{\bf 1.} An obvious consequence is Gromov's non-sqeezing theorem  in any dimension.

{\bf 2.} In the case where the complex structure $J_A$ is integrable Theorem \ref{theo1} can be used as a version of "Hartogs lemma" in order to study the holomorphic convexity properties. 

{\bf 3.} Gromov proved the above result in a slightly different form. For a fixed torus he establishes an existence of a non-constant disc whose boundary contains  a given point of this torus. This does not give an information about the behavior of interior points of the discs. Therefore, it does not allow to determine what subset of $\C^n$ is swept by discs when one varies the boundary tori. For this  reason  we present here a modified version.

\subsection{Manifolds of discs, Fredholm maps and the Sard -Smale theorem}

In this section we outline the proof. 

Fix  a smooth map $\gamma$ from the unit circle to the space of matrix functions satisfying the above assumptions (i), (ii)  such that $\gamma(e^{i0}) = 0$ and $\gamma(e^{i\pi}) = A$. The image by $\gamma$ of the unit circle  is denoted by  ${\cal M}$. Without loss of generality assume  $R = 1$. Denote by $\Lambda^t$ the torus 
$$\Lambda^{t} = \{ (z,w) \in \C^n : \vert z \vert =1, \vert w_j \vert^2 = t_j,
j = 2,...,n \} = b\D \times t_2^{1/2} b\D \times ... \times t_{n}^{1/2} b\D$$
Since $\Lambda^t$ is Lagrangian for $\omega$ and $J$ is tamed by $\omega$, this torus is totally real with respect to $J$.

Denote by $X$ the set $(Z,t)$ of maps $Z: \D \to \C^{n}$ of class $C^{k+1,\alpha}(\D)$ 
and $t= (t_2,...,t_{n})$, $t_j > 0$
satisfying the following assumption:

\bigskip

(iii) {\it  (Boundary data condition)}  For every $(Z,t)$  the boundary condition (\ref{BdCn})
holds. Geometrically this means that $X$ is formed by the smooth discs with
boundaries  attached to the tori $\Lambda^t$ when $t$ runs over
$(\R_+^*)^{n-1}$.

\bigskip

 Denote by $(Z^0,t^0) \in X$ the map $(z^0,w^0)$ where $w^0= b$ is a constant map 
and $z^0$  is a conformal automorphism of $\D$ satisfying $z^0(0) = a$. Put $t^0 =(t^0_2,...,t^0_{n})$ with $t^0_{j} = \vert b_j \vert^2$. Then $Z^0(b\D) \subset \Lambda^{t^0}$ and the disc $Z^0$ satisfies the Cauchy-Riemann equations (\ref{CR1}) with $A = 0$. 

\bigskip

(iv) Denote by $X_0$ a
subset of $X$  of maps satisfying the {\it normalization condition} 
\begin{eqnarray}
\label{normalization}
Z(0) = (z(0),w(0)) = (a,b),\quad z(1) = 1 
\end{eqnarray}

\bigskip

Consider the subset $Y \subset X_0 \times {\cal M} \times C^{k,\alpha}(\D)$ which consists of all $(Z,t,A,h) \in X_0
\times {\cal M} \times C^{k+1,\alpha}(\D) $ satisfying   on $\D$ the non-homogeneous Cauchy-Riemann equations

\begin{eqnarray}
\label{CR2}
Z_{\bar\zeta} - A(Z)\bar Z_{\bar\zeta} - h = 0 
\end{eqnarray}
with the boundary conditions (\ref{BdCn}). Finally, denote by $Y_0$ a subset of $Y$ formed by $(Z,t,A,h)$   homotopic to  $(Z^0,t^0,0,0)$   throught  $Y$. 
In what follows  we assume that $h$ belongs to an open neigbborhood $\Omega$ of the origin in $C^{k,\alpha}(\D)$ which can be shrunk and assumed to be small enough during the proof.

\bigskip

Recall that a bounded linear map $u: E \mapsto E'$ between two Banach
spaces  is called {\it a Fredholm operator} if  $ker u$ and $coker u$ are  finite
dimensional. The number $ind(u) = \dim ker u - \dim coker
u$ is called {\it  the index} of $u$. It is stable under small perturbations of $u$
and a homotopy in the space of Fredholm operators.

A $C^1$-map $F: M \to M'$ between two 
Banach manifolds is called {\it a Fredholm map } if for every point $q \in M$ the tangent
map $dF_q: T_qM \to T_{F(q)}M'$ is Fredholm. The index of every tangent
map  is called the index of $F$; it is denoted $ind (F)$.

A point $q \in M$ is called {\it regular} if the tangent map $dF(q)$ at this point is surjective. A point $q' \in M'$
is called a {\it regular value} of $F$ if the preimage $F^{-1}(q')$ is empty or consists of regular points.

 Consider now the natural  projection 

$$F: Y_0 \to {\cal M} \times \Omega$$
defined by 
$$F:(Z,t,A,h) \mapsto (A,h)$$

The first technical step  is

\begin{prop}
\label{FredholmPro}
\begin{itemize}
\item[(i)] $Y_0$ is a Banach manifold.
\item[(ii)]  The projection $F: Y_0 \to {\cal M} \times \Omega$
is a  Fredholm map with $ind (F) = 0$.
\end{itemize}
\end{prop}

Since this statement is a variation of well-known results \cite{Al,Aud,Gr} which follow from the classical theory of linear integral equations \cite{MiPr,Wend}, we drop the proof. For reader's convenience, we include in  Appendix a proof for the model case of the standard complex structure  used in the next section. Here we only point out that the index is invariant with respect to a homotopy and it suffices to compute it for the above disc $Z^0$. But it is easy to check that  $ind(dF(Z^0)) = 0$.

 The main step of the proof of Theorem is the following

\begin{prop}
\label{proper0}
There exists an open neighborhood of the origin $\Omega$ in $C^{k,\alpha}(\D)$
such that the restriction $F: Y_0 \cap F^{-1}({\cal M} \times \Omega) \to {\cal M} \times \Omega$ is a proper map.
\end{prop}

Admitting for a moment Proposition \ref{proper0}, we  prove Theorem \ref{theo1}.
The key ingredient is provided by  the following general topological principle  due to Smale \cite{Sm}.

\begin{prop} 
(Sard-Smale's theorem.) Let $\tau: M_1 \to M_2$ be a proper Fredholm map between two Banach manifolds. 
Then the set of its regular values is dense in $M_2$. 
For every regular value $p \in M_2$ the preimage $\tau^{-1}(p)$ is a manifold of dimension equal 
to the Fredholm index of $d\tau_q$ (or empty), $q \in \tau^{-1}(p)$. 
Furthermore, for any two regular values $p_1$ and $p_2$ the manifolds $\tau^{-1}(p_1)$ and $\tau^{-1}(p_1)$
 are (non-orientedly) cobordant. 
\end{prop}

The cobordance here means that the union $\tau^{-1}(p_1) \cup \tau^{-1}(p_2)$ is the (non-oriented) boundary $\partial N$ (here we prefer to use the homological notation) of a submanifold $N \subset M_1$.

{\bf Proof of Theorem \ref{theo1}:} The point $(0,0) \in {\cal M} \times \Omega$ is a regular value of $F$ and the preimage $F^{-1}(0,0)$ consists of the single point  $\{ (Z^0,t^0,0,0)\}$. Let $(A,h)$ be another regular value of $F$. Since the set $F^{-1}(0,0)$ consists of a single point, it can not be cobordant to the empty set.  Therefore the preimage $F^{-1}(A,h)$ is not empty. 
Hence the image of $F$ contains a dense subset of ${\cal M} \times \Omega$. Since $F$ is proper, we conclude that $F$ is surjective. In particular   $F^{-1}(A,0)$ is not empty. $\blacksquare$

The remainder of the section is devoted to the proof of Proposition \ref{proper0}.

\subsection{Structure lift and symplectic area } When $A$ is prescribed, the solutions to the non-homogeneous equations (\ref{CR2}) can be viewed as complex discs for a suitable structure $J_{A,h}$ determined by  
$A$ and $h$ ( we will drop $A$ and write just $J_h$ when   $A$ is
 fixed).  The structure $J_h$ is defined on $\D \times \C^n \subset \C^{n+1} = \C_{z_0} \times \C^n_Z$. Setting $z_0(\zeta) = \zeta$, we see that (\ref{CR2}) can be written in the form

\begin{eqnarray}
\label{CR3}
\left\{
\begin{array}{cccc}
& (z_0)_{\bar\zeta} = 0,\\
& Z_{\bar\zeta} - A(Z)\bar Z_{\bar\zeta} - h(z_0)(\bar{z}_0)_{\bar\zeta} = 0
\end{array}
\right.
\end{eqnarray}
We view this PDE system as   the Cauchy-Riemann  equations for a $J_h$-complex disc $\hat Z: \zeta \mapsto (z_0(\zeta),Z(\zeta)) = (\zeta,Z(\zeta))$. This defines the complex matrix of some almost complex structure  which we denote by $J_h$. Hence  $Z = Z(\zeta)$ is a solution of (\ref{CR2}) if and only if $\hat Z$ is a $J_h$-complex disc.

If $Z \in Y$ then its lift $\hat Z$ is glued to the torus $\hat \Lambda^t: = b\D \times \Lambda^t$. The symplectic form $\omega$ lifts to $\C^{n+1}$ as $\hat \omega = (1/2i)d z_0 \wedge d\bar{z_0} + \omega$ i.e. as a standard symplectic form on $\C^{n+1}$. The  torus $\hat \Lambda^t$
remains Lagrangian with respect to $\hat\omega$ and totally real with respect
to $J_h$. We note (see below) that we will consider only $h$ which are close
to the zero-function in the $C^{k,\alpha}$-norm, so the almost complex structure $J_h$
remains tamed by $\hat\omega$.

\bigskip

We use the notation $z = x + iy$ and $w_j = u_j + i v_j$.  Set $\lambda_1=  (1/2)(xdy - ydx)$ and  $\lambda_j = (1/2)(u_jdv_j - v_jdu_j)$. Finally put 
$\lambda = \sum_{j=1}^n \lambda_j$. Then $\omega = d\lambda$. Let $Z \in Y_0$ be a  disc homotopic to $Z^0$. By Stokes' formula  $$area(Z) = \int_{b\D} Z^*\lambda = \pi$$ because the winding numbers of the complex functions  $z$ and $w_j$, $j=2,...,n$ are equal to $1$ and $0$ respectively. In particular, we have the following
\begin{prop}
\label{ProArea1}
The $\omega$-area of every  disc $Z \in Y_0$  is equal to $\pi$.
\end{prop}

 The area of every lift $\hat
Z$ is equal to $2\pi$ since an additional  integral of the $z_0$-component arises. We obtain the following

\begin{prop}
\label{ProArea2}
 Let $(Z,A,h) \in Y_0$. Then the $\hat\omega$ area of $\hat Z$ is equal to $2\pi$.
\end{prop}

 The first consequence is 

\begin{lemma}
Let $(A^m,h^m)$ be a sequence converging in ${\cal M} \times
\Omega$ and  let $(Z^m,t^m,A^m,h^m)$ be  in $Y_0 \cap F^{-1}(A^m,h^m)$. The sequence $(t^m)$ is bounded.
\end{lemma}

\proof  Suppose by absurd that after extracting a subsequence we have $t_m \to \infty$. Then the  lifts $\hat Z^m$ have a bounded area
and unbounded diameters because $\hat Z^m(0) = (0,a,b)$; by the diameter we mean here the maximum on the disc of the distance to the boundary of this disc. But this is impossible, since   the taming assumption (i) implies an upper  bound on the diameter  of every disc in terms of  its area, see \cite{Aud}   $\blacksquare$

\subsection{Separation }  The key statement is 
 
\begin{prop}
\label{separation}
 Fix $\eta > 0$. There exists $\varepsilon_0 > 0$ and a neighborhood $\Omega$ of the origin in $C^{k,\alpha}(\D)$ such that
 for all $(Z,t,A,h)  \in Y_0 \cap F^{-1}({\cal  M} \times \Omega)$ with   $t_j \geq \eta$, $j=2,...,n$ the following estimate holds
$$\vert w_j(\zeta) \vert \geq \varepsilon_0,  j=2,...,n$$ 
for every $\zeta \in \overline\D$. 
\end{prop}

\proof The assertion holds for $A = 0$ and $h = 0$. Using the homotopy assumption in the definition of $Y_0$ and arguing by absurd, assume that there exists a sequence $(Z^m,t^m,A^m,h^m)\in Y_0 \cap F^{-1}({\cal M} \times \Omega)$ such that $t_j  \geq \eta$, $j= 2,...,n$, $A^m \to A^\infty$ (recall that the loop ${\cal M}$ is compact) and $h^m \to 0$, but for some $j$ one has $0 < \inf_\D \vert w_j^m \vert$ for all $m$ and  $\inf_\D \vert w_j^m \vert \to 0$ as $m \to \infty$. Recall that the standard symplectic form $\omega$ is exact and spherical bubbles can not arise. By Gromov's compactness theorem ( extracting a subsequence)
the sequence of lifts $\hat Z^m(\D)$ converges in the Hausdorff distance to a finite union of $J_{A^\infty,0}$-complex discs of class $C^{k,\alpha}(\bar\D)$ with boundaries glued to the torus  $\Lambda^{t^\infty}$, $t^\infty = \lim_{m \to \infty} t^m$. Given the limit disc, after a suitable reprametrization by a sequence of conformal automorphisms, the convergence is in every $C^l(Q)$-norm on every compact subset $Q$ of $\bar \D \setminus \Sigma$, where $\Sigma$ is at most a finite subset of $b\D$. Then the projection of one of the limit discs on $\C^n(z,w)$, say, $ Z$, touches the $J_{st}$-complex hyperplane $P = \{ (z,w): w_j = 0 \}$ at some point $q$.

{\it Case (A).} $q$ is a boundary point of $Z$. Then $t^\infty_j = 0$. Let $\hat Z$ be the "principal" limit disc i.e. the sequence $\hat Z^m$ converges to this disc without additional reparametrization off a finite set.  Outside this finite set, the boundary of the limit disc $Z^\infty$ (the projection of $\hat Z$ to $\C^n$) coincides with the circle $\{ \vert z \vert = 1 \} \times \{ 0 \}$.  Since the almost complex structure is standard near $P$, by the boundary uniqueness theorem for holomorphic functions,  an open subset of $Z^\infty$ is contained in the $J_{A^\infty,0}$-complex hyperplane $P$. Then this inclusion holds globally: a contradiction to the normalization condition (\ref{normalization}).

{\it Case (B).}  $q$ is an interior point.Then we can assume by Case (A) that $t_j^\infty > 0$. Consider on $\D_{z_0} \times \C_{w_j} \subset \C^2$ the equations

\begin{eqnarray}
\label{CR2012}
 \left\{ \begin{array}{ll}
(z_0)_{\overline\zeta} = 0,\\
(w_j)_{\overline\zeta} - h^m(z_0)(\overline z_0)_{\overline\zeta}= 0,
\end{array} \right.
\end{eqnarray}

This is the Cauchy-Riemann equations (\ref{holomorphy}) associated to  an almost complex structure $J^m$ on $\D \times \C$. The  complex matrix of $J^m$ is

\begin{eqnarray}
\left(
\begin{array}{cll}
0 &  0\\
h^m & 0 
\end{array}
\right)
\end{eqnarray}
 
The sequence $\hat Z^m$ converges to $\hat Z$, $q = Z(\zeta_0)$, $\hat q = (\zeta_0,q) = \hat Z(\zeta_0)$. Since $A^\infty = 0$ in a neighborhood of the hyperplane $P$,  the projections $\hat Z^m_j:= (z_0^m,w_j^m)$ satisfy the equations (\ref{CR2012})  near the point $(\zeta_0,0) \in \C^2$ ,i.e. they are $J^m$-complex curves there. Fix $m$ big enough. By the Nijenhuis-Woolf theorem a fixed neighborhood of the disc $\D \times \{ 0 \}$ in $\C^2$ is foliated by a complex 1-parameter family of  $J^m$-complex discs (small deformation of the family $w_j = const$ converging to this family when $m \to \infty$). Then $\hat Z^m_j$ touches one of these discs: a contradiction to the positivity of intersections. $\blacksquare$

\subsection{Proof of Proposition \ref{proper0}}

Now we proceed quite similarly to the case of dimension 2. For reader's convenience we include details.
 
 Arguing by absurd, assume that there exists a sequence $\Omega_j$ of open neighborhoods of the origin converging to the origin, such that for every $j$ the map  $F: Y_0 \cap F^{-1}({\cal M} \times \Omega_j) \to {\cal M} \times \Omega_j$ is not proper. Then for every $j$, there exists a sequence $(A^{m,j},h^{m,j})_m$  converging in ${\cal M} \times
\Omega_j$ to some $(A^{\infty,j},h^{\infty,j})$ as $m \to \infty$ and sequence $(Z^{m,j},t^{m,j},A^{m,j},h^{m,j})_m$ in $Y_0 \cap F^{-1}(A^{m,j},h^{m,j})$ which does not admit a converging subsequence.
Therefore,  by Gromov's compactness theorem \cite{Gr}, a boundary 
disc-bubble arises in every  sequence $(\hat Z^{m,j})_m$ (recall that there are no spherical bubbles since $\omega$ is exact). One can assume that  for every $j$ the sequence $(\hat Z^{m,j}(\D))_m$ converges in the Hausdorff distance to a connected finite union of $J_{A^{\infty,j},h^{\infty,j}}$-complex discs with boundaries glued to the torus $\hat \Lambda^{t^{\infty,j}}$ where $t^{\infty,j}$ is the limit of $(t^{m,j})$.

Choose some $j$. We have  the "principal" limit disc
$\hat Z^{\infty,j}$ which is the limit of the sequence of maps $(\hat Z^{m,j})_m$  converging  as $m \to \infty$ on $\bar\D$ off at most a finite number of boundary points where disc-bubbles arise.   The disc $Z^{\infty,j}$  is centered at $(a,b)$, its
boundary is attached to the torus $\Lambda^{t^{\infty,j}}$. Since the $w$- components of all discs of our sequence are uniformly
separated from the origin  by Proposition \ref{separation},   the limit disc has the same property. Furthermore, its $z_0$-component remains equal to $\zeta$. Its $z$-component is a smooth function on $\bar\D$,  $z:b\D \to b\D$ and  $z(0) = a$. 

\begin{lemma} The winding number of the $z$-component of $\hat Z^{\infty,j}$ does not vanish for all $j$ big enough.
\end{lemma}
\proof Suppose by absurd that it does. Then the area of every disc $\hat Z^{\infty,j}$ is equal to $\pi$ and we  apply Gromov's compactness theorem to the sequence $(\hat Z^{\infty,j})_j$. Since ${\cal M}$ is compact, one can assume that $A^{\infty,j}$ converges to $A^\infty$; we also can assume that $(h^{\infty,j})_j$ converges to $0$ and $(t^{\infty,j})_j$ converges to $t^\infty$. Again we consider the "principal" limit disc $\hat Z^\infty(\zeta) = (\zeta,Z^\infty(\zeta))$ being the limit of the sequence of maps $(\hat Z^{\infty,j})_j$ off at most a finite subset in the boundary.
Then $area(\hat Z^\infty) = \pi$. By the compactness theorem  the sum of this area and the areas of eventually arising  bubbles is equal to $\pi$ and every  bubble has a non-zero area. Hence the bubbles do not arise. Therefore, passing to a subsequence, we have a convergence in 
the $C^{k+1,\alpha}$ norm on the closed disc. Then the disc $Z^\infty$ is $J_{A^\infty}$-complex (because $h^\infty = 0$), its boundary is glued to $\Lambda^{t^\infty}$ and its area is equal to zero.  So it is the constant map equal to some point of $\Lambda^{t^\infty}$. This is a contradiction, since its $z$-component still satisfies $z(0) = a$. $\blacksquare$

Going back to the disc $\hat Z^{\infty,j}$, we conclude that   the winding number of its  $z$-component  is  strictly positive. Indeed, it can not be negative since the symplectic area of the non-constant $J_{A^{\infty,j},h^{\infty,j}}$-complex disc $\hat Z^{\infty,j}$ must be positive. This implies that  
$area(\hat Z^{\infty,j}) = 2\pi$.  On the other hand,  once again by Gromov's compactness
theorem the sum of $area(\hat Z^{\infty,j})$ and the areas of all bubbles  is equal to $area (\hat Z^{m,j}) = 2\pi$. We see that the area of every bubble which could arise in the sequence
$(\hat Z^{m,j})_m$ must be 
equal to $0$  implying that this bubble is constant. This is  a contradiction since bubbles can not be constant.  $\blacksquare$

\section{Gluing a complex  disc to a Lagrangian submanifold of $\C^n$}

An adaptation of the method employed in previous section allows to prove the following  result \cite{Gr}:

\begin{thm}
\label{GluingDiscsTheo1}
Let $E$ be a smooth compact Lagrangian submanifold in $(\C^n,\omega_{st})$. Then there exists a non-constant $J_{st}$-holomorphic disc attached to $E$.
\end{thm}

Our exposition here is inspired by \cite{Al,IvSh}. As above  the standard symplectic form $\omega_{st}$ is denoted simply by $\omega$. In the following we only slightly modify the notation of previous section.

\subsection{Adapted manifolds of discs.}  Fix a point $a \in E$ and consider the constant 
holomorphic map $f^0(\zeta) \equiv a$. Fix $k \geq 1$ and $0 < \alpha < 0$.
 Consider the set  $X \subset C^{k,\alpha}(\D)$ 
consisting of maps $f:\D \to \C^n$  with $f(b\D) \subset E$ and $f(1)= a$.
 Assume in addition that $E = \{ \rho = 0 \}$ where $\rho = (\rho_1,...,\rho_n):\C^n \to \R^n$ 
is a smooth map of maximal rank. Furthermore 
\begin{eqnarray}
\label{TotReal}
\partial \rho_1 \wedge ... \wedge \partial \rho_n \neq 0
\end{eqnarray}
because $E$ is totally real. Then a holomorphic map $f \in C^{k,\alpha}(\D)$ is in $X$ if and only if it 
is a solution of the following non-linear Riemann-Hilbert type boundary value problem: 

\begin{displaymath}
(RH): \left\{ \begin{array}{ll}
\frac{\partial f(\zeta)}{\partial\overline\zeta} = 0, \,\,\,\zeta \in \D\\
\rho(f)(\zeta) = 0, \,\,\, \zeta \in b\D\\
f(1)= a
\end{array} \right.
\end{displaymath}

Denote by $X_0$ the subset of  $X$ formed by the discs homotopic in $X$ to the constant disc $f^0$. Set $G = C^{k,\alpha}(\D)$. Consider  the cartesian product $X_0 \times G$ and define a subset 
$$Y \subset X_0 \times G = \{ (f,h) : \frac{\partial f}{\partial \overline\zeta} = h \}$$
 Denote by $F: X_0 \times G \to G$ the natural projection.

As in the previous section we have

\begin{prop}
\begin{itemize}
\item[(i)] $Y$ is a Banach manifold.
\item[(ii)]  The projection $F: Y \to G$
is a  Fredholm map with $ind (F) = 0$.
\end{itemize}
\end{prop}
\proof Since $X_0$ consists of discs homotopic to the constant map $f^0 \equiv a \in E$, 
the Fredholm index of $dF_{(f,h)}$ is independent of $(f,h)$ and coincides with the index 
of $dF_{(f^0,0)}$ which is equal to $0$. $\blacksquare$

A crucial property of $F$ is given by the following 

\begin{lemma}
\label{non-sur}
 The projection $F:Y \to G$ is not surjective. 
\end{lemma}
\proof Suppose by contradiction that $F$ is surjective. 
Then for every $t > 0$ there exists $f^t \in X_0$ such that $(f^t,h^t) \in Y$ where $h^t(\zeta):= (t,0,...,0)$.
 On the other hand, $\partial f^t/\partial\overline \zeta = h^t$ and, in particular $\partial f^t_1/\partial\overline\zeta = t$. 
Hence $f^t_1 = t \overline \zeta + q^t(\zeta)$ where $q^t$ is a function holomorphic in $\D$. 
Since $f^t(b\D) \subset E$, the family $(f^t)$ is bounded on $b\D$ by a constant $C > 0$ 
independent of $t$. Therefore $\vert \overline\zeta + t^{-1}q^t(\zeta) \vert \leq t^{-1}C$ 
for $\zeta \in b\D$. However the function $\zeta \mapsto \overline\zeta + t^{-1}q^t(\zeta)$ 
is harmonic in $\D$ and by the maximum principle a similar estimate holds for all $\zeta \in \D$. Letting $t \to \infty$, we obtain that the function $\zeta \mapsto \overline\zeta$ can be uniformly approximated by holomorphic 
functions in $\D$: a contradiction. $\blacksquare$

\subsection{Non-linear Fredholm alternative.} In order to conclude the proof of Theorem \ref{GluingDiscsTheo1}, it suffices to establish the following

\begin{prop}
\label{FredholmAlter}
Suppose that there does not exist a non-constant Bishop disc for $E$. Then the projection $F: Y \to G$ is surjective.
\end{prop}
Then  the theorem follows by contradiction. The remainder is devoted to the proof of Proposition.

We deal with the non-homogeneous Cauchy-Riemann equations

\begin{eqnarray}
\label{d-bar}
\frac{ \partial f}{\partial\overline\zeta}(\zeta) = h
\end{eqnarray}
on the unit disc.In order to use Gromov's compactness theorem, we again view its solutions as $J_h$-complex discs
for a suitably choosen almost complex structure $J_h$ in $\C^{n+1}$. Quite similarly to previous section (cf. with (\ref{CR2012})) , given $h$   
consider the almost complex structure $J_h$ on $\D \times \C^n \subset \C^{n+1}$. The equations (\ref{d-bar}) are equivalent to the fact that the lift
 $\hat f : \zeta \mapsto (\zeta,f(\zeta))$ of $f$ is a $J_h$-complex disc. 

\begin{lemma}
\label{proper}
The projection $F: Y \to G$ is proper. 
\end{lemma}
\proof We must show that for every sequence $h^k \to h^\infty$ in $G$
 and every sequence $(f^k)$ such that  $(f^k,h^k) \in Y$ there exists a subsequence of $(f^k)$ converging in $X_0$.
The structures $J_{h^k}$ converge to $J_h$ and the 
discs $\hat f^k(\zeta) = (\zeta,f^k(\zeta))$ are $J_{h^k}$-complex. 
Their boundaries are attached to the manifold $\hat E:=b\D \times E$. 
Since $E$ is a Lagrangian manifold, it follows that  $\hat E$ is a Lagrangian manifold in $\C \times \C^n$ with respect to the 
symplectic form $\hat \omega = C\frac{i}{2}dz_0 \wedge d\overline z_0 + \omega$.
 The sequence $(h^k)$ is bounded, which implies that there exists a constant $C > 0$ 
such that  the structures $J_{h^k}$ is tamed by $\hat \omega$ for all $k = 0,1,...,\infty$.

Denote by $\lambda$ a primitive of $\hat\omega$. Then
$$\int_{\hat f^k(\D)} \hat \omega = \int_{\hat f^k(b\D)}  \lambda$$
On the other hand the boundaries $\hat f^k(b\D)$ of our discs are homotopic 
(the discs are homotopic) and $d\lambda\vert_{\hat E} = \hat\omega \vert_{\hat E} = 0$. 
Then by Stokes' theorem the last integral is independent of $k$. 
Gromov's  compactness theorem implies that there are only the folowing possibilities: 

(a)  The limit of some subsequence of $(\hat f^k)$ contains a disc-bubble $\psi$. Since $\hat f^k$ are the graphs over $\D$, it follows easily from the definition of disc-bubbles that 
$\psi$ is ``vertical'', i.e. has the form $\psi:\zeta \mapsto (q, f(\zeta))$ where $q$ is a point of the unit circle $b\D$. Indeed, one readily sees from the definition of a disc-bubble that the renormalizing sequence $\phi^n$ of conformal biholomorphisms is not compact (in the compact-open topology) and converges to a boundary point of the unit disc. Then it follows from the Cauchy-Riemann equations associated with $J_{h^\infty}$ that $f$ is a usual  holomorphic (with respect to $J_{st}$) disc 
 attached to $E$ (cf. with the equations (\ref{CR2012}) of previous section). The disc $f$ is non-constant because  the  bubble $\psi$ is non-constant. This contradicts  the assumption of Proposition \ref{FredholmAlter}. 

(b) The limit of some subsequence of $({\hat f}^k)$ contains a non-constant $J_{h^\infty}$-complex sphere. This is impossible since $\hat\omega$ is exact.

Thus, only the last possibility realizes:

(c)  there exists a subsequence converging in $C^{k+1,\alpha}(\D)$-norm. $\blacksquare$.

Now using the Sard-Smale theorem we conclude as in previous section that 
 $F$ is surjective. This proves Proposition \ref{FredholmAlter} and Theorem 
\ref{GluingDiscsTheo1}. $\blacksquare$

\subsection{Exotic symplectic structures} As a consequence we obtain the existence of exotic symplectic structures on $\R^{2n}$. Denote by $\omega_{st}$ the standard symplectic form of $\C^n$, identyfing $\C^n$ with $\R^{2n}$. Set $\lambda_{st} = \sum_j x_jdy_j$ so that $d\lambda_{st} = \omega_{st}$.

 A symplectic structure $\omega$ on $\R^{2n}$ is called exotic if there is no  global diffeomorphic map $\phi:\R^{2n} \to \R^{2n}$
such that $\phi^*\omega_{st} = \omega$.

\begin{cor}
Let $E$ be a compact Lagrangian submanifold in $(\C^n, \omega_{st})$. Then the restriction $\lambda_{st}\vert_E$ represents a non-zero class in $H^1(E,\R)$.
\end{cor}
\proof Since $d\lambda_{st}\vert_E = 0$, then $\lambda_{st}\vert_E$ does represent a class 
in $H^1(E,\R)$. If this class is zero, then for every closed smooth curve $\gamma \subset E$ 
we have $\int_\gamma \lambda_{st} = 0$. 

Let $f$ be a non-constant Bishop disc glued to $E$ and let $\gamma = f(b\D)$. Then by Stokes' formula
$$ \int_{b\D} f^*\lambda_{st}  = \int_\D f^*\omega_{st} = area[f(\D)] > 0$$
which proves the corollary. $\blacksquare$

How to find an exotic structure? Consider the standard torus $\Lambda = \{ z = (z_1,z_2) \in \C^2: \vert z_j \vert = 1  \}$ which is Lagrangian for $\omega_{st}$ in $\C^2$.

\begin{lemma}
Suppose that $\omega$ is a symplectic form such that for some $1$-form $\lambda$ we have $\omega = d\lambda$ on $\R^4$ and $\lambda \vert_\Lambda = 0$. Then $\omega$ is exotic.
\end{lemma}
\proof Suppose that $\phi$ is a diffeomorphism of $\R^4$ satisfying $\omega = \phi^*\omega_{st}$. Then $\phi^*\lambda_{st} - \lambda$ is a closed 1-form and there exists a function $h$ on $\R^4$ such that $dh = \phi^*\lambda_{st} - \lambda$. Then 
$$\phi^*\lambda_{st}\vert_\Lambda = (\phi^*\lambda_{st} - \lambda) \vert_\Lambda = d(h\vert_\Lambda)$$
so that $\phi^*\lambda_{st}\vert_\Lambda$ is exact. Therefore $\lambda_{st}\vert_{\phi(\Lambda)}$ is exact and $\phi(\Lambda)$ is Lagrangian for $\omega_{st}$: a contradiction. $\blacksquare$.

It turns out that it is not difficult to write explicitely a symplectic structure satisfying the assumptions of the above Lemma, see \cite{MS0}.

\section{Appendix: Fredholm property}
Using notations of previous section, we  prove here that $F:Y \to G$ is a Fredholm map.Let $(f_0,h_0) \in Y$.  We follow \cite{Al}. The tangent space $T_{(f_0,h_0)}Y$ to $Y$ at $(f_0,h_0)$ is formed by the maps $(\dot f,\dot h) \in C^{k+1,\alpha}(\D) \times C^{k,\alpha}(\D)$ satisfying
\begin{displaymath}
 \left\{ \begin{array}{ll}
\frac{\partial \dot f}{\partial \overline\zeta} = 0, \zeta \in \D,\\
2\Re P(\zeta)\dot f(\zeta) = 0, \zeta \in b\D,\\
\dot f(1) = 0
\end{array} \right.
\end{displaymath}

Here $P(\zeta)$ is the Jacobian matrix 
$$ \left ( \frac{\partial \rho}{\partial Z}(f_0(\zeta)) \right )$$
Since $\dot h$ is arbitrary, we identify $T_{(f_0,h_0)}H$ with the space of maps $\dot f:\D \to \C^n$ satisfying $\dot f(\zeta) \in T_{f_0(\zeta)}(E)$ for $\zeta \in b\D$ and $\dot f(1) = 0$. We identify $G$ with the tangent space $T_{h_0}G$. Then the tangent map $dF: T_{(f_0,h_0)}H \to G$ to $F$ at $(f_0,h_0)$
is 
$$dF: \dot f \mapsto \dot h = \frac{\partial \dot f}{\partial \overline\zeta}$$
The condition that $E$ is totally real is equivalent to the fact that $\det P(\zeta) \neq 0$, $\zeta \in b\D$ ("the Lopatinski condition"). Set $A = C^{k+1,\alpha}(\D,\C^n)$ and $B = C^{k,\alpha}(\D,\C^n) \times C^{k+1,\alpha}(b\D,\C^n)$. Consider the linear operator
$$L:A \to B,$$
$$L: \dot f \mapsto (\frac{\partial \dot f}{\partial \overline\zeta}, \Re P\dot f\vert_{b \D})$$
According to   \cite{Wend}, Th. 3.2.5., the operator $L$ is Fredholm.Define the operator 
$L_1: A \to B \times \C$ by $L_1(\dot f) = (L(\dot f),\dot f(1))$. Obviously $L_1$ is Fredholm. Let $I :T_{(f_0,h_0)}Y \to A$ be the inclusion map. Consider the map $I':G \to B \times \C$ defined by $I':\dot h \mapsto (\dot h,0,0)$. Then we have the commutative diagramm
$$I' \circ dF = L_1 \circ I$$
In particular, $I(ker (dF)) \subset ker L_1$ and so $ker dF$ has a finite dimension.
Furthermore, $I'(im (dF)) \subset im L_1$ and so the induced quotient map $\widehat {I'}:
G/im (dF) \to (B \times \C)/im L_1$ is correctly defined. Obviously it is injective. Hence $\dim coker dF \leq \dim coker L_1 < + \infty$. This means that $dF$ is a Fredholm map.

\end{document}